\documentclass[submission]{arxivstyle}
\newtheorem{thm}{Theorem}
\newtheorem{lem}{Lemma}
\newtheorem{proposition}{Proposition}
\renewcommand{\geq}{\geqslant}
\renewcommand{\leq}{\leqslant}
\usepackage{tikz}

\DeclareMathOperator{\Jac}{Jac}
\DeclareMathOperator{\Ima}{Im}
\newcommand{\Om}{\Omega}
\newcommand{\ep}{\varepsilon}
\newcommand{\de}{\delta}

\DeclareMathOperator{\Id}{Id}

\title[Finiteness of the group for weighted walks in orthants]{On the finiteness of the group associated with weighted walks in multidimensional orthants}

\author[A. Elvey Price, E. Humbert and K. Raschel]{Andrew Elvey Price\thanks{\href{mailto:andrew.elvey@univ-tours.fr}{andrew.elvey@univ-tours.fr}.  Partially supported by ANR CartesEtPlus (ANR-23-CE48-0018) and ANR
IsOMA (ANR-21-CE48-0007)}\addressmark{1} \and Emmanuel Humbert\thanks{\href{mailto:emmanuel.humbert@univ-tours.fr}{emmanuel.humbert@univ-tours.fr}. Partially supported by ANR Einstein-PPF (ANR-23-CE40-0010).}\addressmark{1} \and Kilian Raschel\thanks{\href{mailto:raschel@math.cnrs.fr}{raschel@math.cnrs.fr}. Partially supported by ANR RAWABRANCH (ANR-23-CE40-0008).}\addressmark{1}}

\address{\addressmark{1}Université de Tours, CNRS, Institut Denis Poisson}

\received{\today}

\abstract{In the study of walks with small steps confined to multidimensional orthants, a certain group of transformations plays a central role. In particular, several techniques to potentially compute the generating function, including the orbit sum method, can only be applied when this group is finite. In this note, we present three new results concerning this group. First, in two dimensions, we provide a complete characterization of the weight parameters that yield a finite group. In higher dimensions, we show that whenever the group is finite, it must necessarily be isomorphic to a simpler reflection group. Finally, in dimension three, we give a full classification of the parameters leading to a finite group that also satisfies an additional Weyl property.}

\resume{Dans l'étude des marches restreintes à des orthants multidimensionnels, il s'avère qu'un certain groupe de transformations joue un rôle crucial, en particulier au travers de son éventuelle finitude. En effet, si ce groupe est fini plusieurs méthodes peuvent être appliquées (comme la somme sur l'orbite) pour calculer la fonction génératrice d'intérêt. Dans cette note nous présentons trois nouveaux résultats sur ce groupe. Tout d’abord, en dimension deux, nous donnons une caractérisation complète des poids conduisant à un groupe fini. En dimension supérieure, nous montrons que si le groupe est fini, il est nécessairement isomorphe à un groupe de réflexions plus simple. Enfin, en dimension trois, nous proposons une classification complète des paramètres menant à un groupe fini possédant en outre une propriété de type Weyl.}


\usepackage[backend=bibtex]{biblatex}
\begin{document}

\maketitle
%% note that you DO NOT have to put your abstract here -- it is generated by \maketitle and the \abstract and \resume commands above

\section{Introduction}

A lattice walk is a sequence of points $P_0, P_1,\ldots ,P_n$ of $\mathbb Z^d$, $d\geq 1$. The points $P_0$ and $P_n$ are its starting and end points, respectively, the consecutive differences $P_{i+1}-P_i$ its steps, and $n$ is its length. Given a set $\mathcal S\subset \mathbb Z^d$, called the step-set, a set $C\subset \mathbb Z^d$ called the domain (which in this paper will systematically be the cone $\mathbb R_+^d$, called the $d$-dimensional orthant), and elements $P$ and $Q$ of $C$, we are interested in the number 
$e_C(P,Q; n)$
of (possibly weighted)\ walks of length $n$ that start at $P=P_0$, have all their steps in $\mathcal S$, have all their points in $C$, and end at $Q=P_n$. Normalising the weights (which are assumed to be non-negative)\ with the condition that they sum to one, we obtain transition probabilities, and the number $e_C(P,Q; n)$ can be interpreted as the probability that a random walk starting at $P$ will reach the point $Q$ at time $n$ while remaining in the domain $C$.

In the last twenty years, there has been a dense research activity in the mathematical community on the enumerative aspects of walks confined to cones, in particular to the $d$-dimensional orthant. To summarise, three main questions have attracted most attention: the first is to determine, if possible, a closed-form formula for the number of walks $e_C(P,Q; n)$. Of course, such an explicit formula is not expected to exist in general, and in most cases can be explained by bijections with other combinatorial objects. The second question concerns the asymptotic behaviour, e.g.\ of the number of walks $e_C(P,Q; n)$, in the regime where the length $n\to\infty$. The last question focuses on the complexity of generating functions associated with these models, such as the series
\begin{equation}
\label{eq:excursions_series}
   \sum_{n\geq0} e_C(P,Q;n)t^n.
\end{equation}
One then asks whether the series satisfies any algebraic or differential equation. Answering this question allows us to classify the models according to the complexity of their generating function. The three questions are by no means independent: for example the possible asymptotic forms of the numbers $e_C(P,Q;n)$ depend heavily on the complexity of the associated generating function.

To present the main results of this paper, we recall a tool that has played a crucial role in obtaining several previous results in the literature: the {\em group} of the walk model. This is a group that was first introduced in two dimensions, in a combinatorial context, in \cite{BMMi-10}, following the idea of Fayolle, Iasnogorodski and Malyshev in \cite{FaIaMa-17}. This group will be properly defined later in the paper, but can be presented informally as a symmetry group of involutions defined by the step-set $\mathcal S$ of the model. The main application is that, if the group is finite, its action on a functional equation naturally associated with the model can lead to explicit expressions for the generating functions \eqref{eq:excursions_series} as well as information on the asymptotics and algebraic complexity. 
%\aep{I think we should be more precise in the following sentences: The method works completely in dimension 2, whereas it only sometimes works in higher dimension, so it doesn't make sense to say this as if it works equally well in dimension 2 and higher.} 
In principle, this method works in dimension two \cite{FaIaMa-17,BMMi-10} and, in higher dimensions, may apply in a few favorable cases \cite{BoBMKaMe-16,MeMi-16,Ya-17,BuHoKa-21}. However, there is no known criterion to determine whether the group is finite (even in two dimensions). This is precisely the question we address in this paper.

\section{Main results}

We work under the small step-set assumption, that is, $\mathcal{S} \subset \{-1, 0, 1\}^d$.
Define the inventory of the model $\mathcal S \in \{-1,0,1\}^d$ as follows:
\begin{equation}
\label{eq:inventory}
    \chi_\mathcal S(x_1,\ldots,x_d) = \sum_{(i_1,\ldots,i_d)\in\mathcal S} w(i_1,\ldots,i_d)x_1^{i_1}\cdots x_d^{i_d},
\end{equation}
where $w(i_1,\ldots,i_d)>0$ is the weight of the step $(i_1,\ldots,i_d)\in\mathcal S$. We will normalize the \textcolor{black}{non-negative} weights in such a way that $\chi_\mathcal S(1,\ldots,1)=1$, so that they are also  transition probabilities. Most of the time we shall assume that the step-set satisfies the following irreducibility assumption:
\begin{enumerate}[label={\rm (H\arabic*)},ref={\rm (H\arabic*)}]
     \item\label{it:hypothesis_irreducible}The step-set $\mathcal S$ is not included in any half-space $\{y\in\mathbb R^d : \langle x,y\rangle\geq0\}$ with $x\in\mathbb R^d\setminus\{0\}$, $\langle \cdot,\cdot\rangle$ denoting the classical Euclidean inner product.
     %For any two points $P,Q$ in the domain $C$, the set $\{n\in\mathbb N : e_C(P,Q;n)\neq 0\}$ is non-empty, with $\mathbb N=\{1,2,\ldots\}$.
\end{enumerate}
%As a consequence of \ref{it:hypothesis_irreducible}, the walk can visit any point in the domain $C$. 
A consequence of \ref{it:hypothesis_irreducible} is that the model is truly $d$-dimensional and is not directed towards a half-space, unlike, for example, the singular walks considered in \cite{BMMi-10}.

\subsection{The combinatorial group}

This group was first introduced in the context of two-dimensional walks \cite{FaIaMa-17,BMMi-10} and turns out to be very useful. Let $\chi_\mathcal S$ be the inventory \eqref{eq:inventory}, and define $A_{j},B_{j},C_{j}$ for $j=1,\ldots,d$ as follows
\begin{align*}
     \chi_\mathcal S(x_1,\ldots,x_d) &=x_1A_{1}(x_2,\ldots,x_d)+B_1(x_2,\ldots,x_d)+\overline{x_1}C_{1}(x_2,\ldots,x_d)\\
                     &=x_2A_{2}(x_1,x_3,\ldots,x_d)+B_2(x_1,x_3,\ldots,x_d)+\overline{x_2}C_{2}(x_1,x_3,\ldots,x_d)\\
                     &=\cdots\\
                     &=x_dA_{d}(x_1,\ldots,x_{d-1})+B_d(x_1,\ldots,x_{d-1})+\overline{x_d}C_{d}(x_1,\ldots,x_{d-1}),
\end{align*}
where $\overline{x_i}=\frac{1}{x_i}$. Under the assumption~\ref{it:hypothesis_irreducible}, the functions $A_{1},\ldots,A_d$ and $C_{1},\ldots,C_d$ are all non-zero. The group of $\mathcal S$ is defined as the group 
\begin{equation}
    \label{eq:group_G_def}
    G=\langle\varphi_1,\ldots,\varphi_d\rangle
\end{equation}
of birational transformations of the variables $[x_1,\ldots,x_d]$ generated by the involutions:
\begin{equation}
\label{eq:expression_generators}
     \left\{\begin{array}{rcl}
     \varphi_1([x_1,\ldots,x_d]) & =& \left[\overline{x_1}\frac{C_{1}(x_2,\ldots,x_d)}{A_{1}(x_2,\ldots,x_d)},x_2,\ldots,x_d\right],\smallskip\\
     \varphi_2([x_1,\ldots,x_d]) & =& \left[x_1,\overline{x_2}\frac{C_{2}(x_1,x_3,\ldots,x_d)}{A_{2}(x_1,x_3,\ldots,x_d)},x_3,\ldots,x_d\right],\smallskip\\
     &\cdots & \\
     \varphi_d([x_1,\ldots,x_d]) & =& \left[x_1,\ldots,x_{d-1},\overline{x_d}\frac{C_{d}(x_1,\ldots,x_{d-1})}{A_{d}(x_1,\ldots,x_{d-1})}\right].
     \end{array}\right.
\end{equation}
The generators of the group $G$ therefore satisfy the relations $\varphi_1^2=\cdots=\varphi_d^2=\Id$, plus possible other relations, depending on the model. Note that if $\mu, \alpha_1, \ldots, \alpha_d > 0$, then the model defined by the inventory $(x_1, \ldots, x_d) \mapsto \mu, \chi_{\mathcal{S}}(\alpha_1 x_1, \ldots, \alpha_d x_d)$ shares the same combinatorial group as the model $\mathcal{S}$. Such transformations are referred to as central weightings.

\subsection{A new complete classification in dimension two} \label{parag_dim_2}
\begin{figure}
\begin{center}
 \begin{tikzpicture}[scale=.7] 
    \draw[->,white] (1,2) -- (0,-2);
    \draw[->,white] (1,-2) -- (0,2);
    \draw[->] (0,0) -- (1,1);
    \draw[->] (0,0) -- (1,0);
    \draw[->] (0,0) -- (-1,-1);
    \draw[->] (0,0) -- (-1,0);
    \draw[->] (0,0) -- (0,1);
    \draw[->] (0,0) -- (-1,1);
    \draw[->] (0,0) -- (1,-1);
    \node at (-1.2,0) {$3$};
    \node at (1.2,1.3) {$15$};
    \node at (-1.2,-1.4) {$2$};
    \node at (0,1.3) {$13$};
    \node at (1.2,0) {$9$};
    \node at (1.2,-1.4) {$6$};
    \node at (-1.2,1.3) {$5$};
   \end{tikzpicture}\qquad
\begin{tikzpicture}[scale=.7] 
    \draw[->,white] (1,2) -- (0,-2);
    \draw[->,white] (1,-2) -- (0,2);
    \draw[->] (0,0) -- (0,1);
    \draw[->] (0,0) -- (1,1);
    \draw[->] (0,0) -- (1,0);
    \draw[->] (0,0) -- (1,-1);
    \draw[->] (0,0) -- (0,-1);
    \draw[->] (0,0) -- (-1,-1);
    \node at (-1.2,-1.4) {$1$};
    \node at (1.2,1.3) {$7$};
    \node at (0,-1.4) {$2$};
    \node at (0,1.3) {$7$};
    \node at (1.2,0) {$5$};
    \node at (1.2,-1.4) {$1$};
   \end{tikzpicture}
   \qquad
\begin{tikzpicture}[scale=.7] 
    \draw[->,white] (1,2) -- (0,-2);
    \draw[->,white] (1,-2) -- (0,2);
    \draw[->] (0,0) -- (1,1);
    \draw[->] (0,0) -- (1,0);
    \draw[->] (0,0) -- (-1,0);
    \draw[->] (0,0) -- (-1,-1);
    \node at (-1.2,0) {$6$};
    \node at (1.2,0) {$2$};
    \node at (1.2,1.3) {$4$};
    \node at (-1.2,-1.4) {$3$};
   \end{tikzpicture}
   \qquad
\begin{tikzpicture}[scale=.7] 
    \draw[->,white] (1,2) -- (0,-2);
    \draw[->,white] (1,-2) -- (0,2);
    \draw[->] (0,0) -- (1,1);
    \draw[->] (0,0) -- (1,0);
    \draw[->] (0,0) -- (0,-1);
    \draw[->] (0,0) -- (-1,0);
    \draw[->] (0,0) -- (0,1);
    \draw[->] (0,0) -- (-1,1);
    \draw[->] (0,0) -- (1,-1);
    \node at (-1.2,0) {$1$};
    \node at (1.2,1.3) {$1$};
    \node at (0,-1.4) {$1$};
    \node at (0,1.3) {$2$};
    \node at (1.2,0) {$2$};
    \node at (1.2,-1.4) {$1$};
    \node at (-1.2,1.3) {$1$};
   \end{tikzpicture}
\end{center}
\caption{Examples of finite group models taken from \cite{KaYa-15}. From left to right, a model with a group of order $4$, $6$, $8$ and $10$. They should be normalised to satisfy the property that $\sum_{(i,j)\in\mathcal S}w(i,j)=1$. The paper \cite{KaYa-15} actually contains infinite families of finite group examples. For example, the third example defines a group model of order $8$ if and only if the associated weights $w(i,j)$ satisfy $w(1,0)w(-1,0)=w(1,1)w(-1,-1)\neq 0$.% In the examples above, the weights are taken to be rational numbers; note that non-rational weights are also allowed.
%\aep{(commented out comment about weights being rational)}
}
\label{fig:step_sets}
\end{figure}
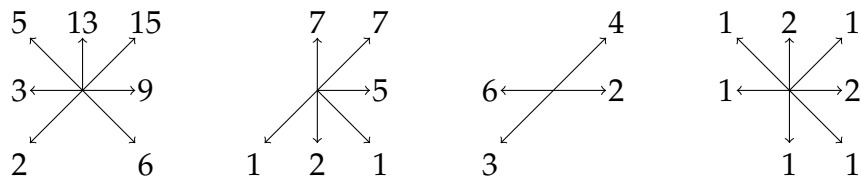

In the unweighted case (that is, when all jumps in the step-set have the same weight), the models were classified in \cite{BMMi-10} according to the finiteness of the group. Among the $79$ non-equivalent models, only $23$ have a finite group. For example, unweighting the first and third models of Figure~\ref{fig:step_sets}, yields models with groups of order $4$ and $8$, respectively.

%First note that to any unweighted model, one can associate a so-called central weighting without changing the group. Central weightings are defined by $c(i,j)=\mu\alpha^i \beta^jw(i,j)$ for some positive $\alpha,\beta,\mu$.
In the general weighted setting, the paper \cite{KaYa-15} characterizes the weights that produce models with groups of order $4$, $6$, and $8$. For example, they find two families of models with a group of order eight: the first one (family 4a in \cite{KaYa-15}) is defined by $w(1,1) = w(0,1) = w(0,-1) = w(-1,-1) = 0$ and $w(1,-1)w(-1,1)= w(1,0)w(-1,0)\neq 0$, the second one is our third example on Figure~\ref{fig:step_sets}.  
%In total, there is one (resp.\ six, two) family of models corresponding to groups of order four (resp.\ six, eight).
Furthermore, the authors identify three isolated models with a group of order $10$, namely the one on the right of Figure~\ref{fig:step_sets}, together with its vertical reflection and its reflection through the origin. Kauers and Yatchak leave open the possibility of additional finite-group models. Essentially, we prove that these are indeed all the finite-group models.

\begin{thm}
\label{thm:1}
In dimension $2$, the only possible orders of $G$ are $4$, $6$, $8$ and $10$.
\end{thm}
In fact, we prove more: first, we recover all the families identified in \cite{KaYa-15}, and we show that the only models with a group of order $10$ are the rightmost model in Figure~\ref{fig:step_sets}, its two symmetric versions, and all their central weightings.

\subsection{A reflection group} \label{reflection_group}
Let $\mathcal S$ be a step-set satisfying~\ref{it:hypothesis_irreducible}, and let $\chi_\mathcal S$ be its inventory~\eqref{eq:inventory}. The system of equations
\begin{equation}
\label{eq:chi_critical_point}
     \frac{\partial \chi_\mathcal S}{\partial x_1}=\cdots = \frac{\partial \chi_\mathcal S}{\partial x_d}=0
\end{equation}
admits a unique solution in $(0,\infty)^d$, denoted by $\boldsymbol{x_0}$. The point $\boldsymbol{x_0}$ has the following interpretation: if we consider the central weighting of the model $\mathcal S$ with the weights $w(i_1,\ldots,i_d)\frac{\boldsymbol{x_0}^{(i_1,\ldots,i_d)}}{\chi_\mathcal S(\boldsymbol{x_0})}$ using the multi-index notation, we get a zero drift model. Define now
the covariance matrix 
\begin{equation}
\label{eq:expression_covariance_matrix}
     \Delta=\bigl(a_{ij}\bigr)_{1\leq i,j\leq d},\quad \text{with } a_{ij}=\frac{\frac{\partial^2 \chi_\mathcal S}{\partial x_i\partial x_j}(\boldsymbol{x_0})}{\sqrt{\frac{\partial^2 \chi_\mathcal S}{\partial x_i^2}(\boldsymbol{x_0})\cdot \frac{\partial^2 \chi_\mathcal S}{\partial x_j^2}(\boldsymbol{x_0})}}.
\end{equation}
Let $\Delta^{-\frac{1}{2}}$ denote the inverse of the symmetric, positive definite square root of the covariance matrix $\Delta$.
This construction is designed to satisfy the following property: the walk with steps $\Delta^{-\frac{1}{2}}\mathcal S$ (together with the above reweighting) has zero drift and identity covariance matrix. In other words, the matrix $\Delta$ is canonical in the sense that it allows one to transform the random walk into a new walk lying in the domain of attraction of a universal Brownian motion with zero drift and identity covariance matrix.
%\aep{Can we give some intuition for this matrix? Or where it has appeared before? E.g., in the case that the weights sum to $1$ and the model has 0 drift, $a_{ij}$ denotes the covariance between the $i$th and $j$th coordinate of the position after a step.}
Consider the $d$-dimensional unbounded polytope
\begin{equation}
\label{eq:def_new_domain}
    T=\Delta^{-\frac{1}{2}}\mathbb R_+^d.
\end{equation}
We denote the canonical basis of $\mathbb R^d$ by $(e_i)_{1\leq i\leq d}$, so the orthant $\mathbb R^d_+$ is bounded by the hyperplanes  
\begin{equation}
\label{eq:def_Gi}
G_i = \text{Span}\{e_1, \dots, e_{i-1}, e_{i+1}, \dots, e_d\},\quad i \in \{1, \dots, d\},
\end{equation}
and thus $T$ in \eqref{eq:def_new_domain} is bounded by the hyperplanes \(H_i = \Delta^{-\frac{1}{2}} G_i \). 
Finally, we introduce the reflection group
\begin{equation}
    \label{eq:group_H_def}
    H = \langle r_1,\ldots,r_d\rangle,
\end{equation}
where $r_{i}$ is the orthogonal reflection in side $H_{i}$ of $T$.
%By definition, this is the group generated by the reflections with respect to the sides of the linear transformation of the orthant \( \Delta^{-\frac{1}{2}} \mathbb R_+^d \). 

In \cite[Cor.~6]{GoHuRa-25}, the authors prove the existence of a surjective morphism that sends the set of generators $(\varphi_1, \ldots, \varphi_d)$ of $G$ in \eqref{eq:expression_generators} to the set of generators $(r_1,\ldots,r_d)$ of $H$. In this note, we present a strong refinement of that previous result: 

\begin{thm}
\label{thm:2}
If $G$ is finite, then $G$ is isomorphic to $H$.
\end{thm}
In particular, this implies that $G$ is a reflection group, which is convenient as reflection groups are well understood.

\subsection{Classification of finite groups in three dimensions and beyond}

 Theorems \ref{thm:1} and \ref{thm:2} are not only interesting in their own right, but can also be combined to completely determine, in arbitrary dimension, the cases where $G$ is finite and the cone $T$ in \eqref{eq:def_new_domain} is a {\em Weyl chamber} of the reflection group $H$. % (for instance, a Weyl chamber of type $B$ is $W_B = \{0 < x_1 < x_2 < \cdots < x_d\}$).
 In these case we will say that $G$ has the Weyl property. Focusing on these models is natural, as we believe that they are precisely the cases where the reflection principle and orbit-sum methods can be used to directly compute the number of walks; see \cite{Fe-14,MeMi-16,Ya-17,BuHoKa-21}.

We do not state a precise general result, but we illustrate how this can be done in a specific case by classifying all irreducible groups in dimension 3 that satisfy the Weyl property.

\begin{thm}
\label{thm:3}
In dimension $d=3$, the models that give rise to finite groups with the Weyl property can be completely classified. See the end of Section~\ref{sec:5} for a concrete example.
\end{thm}

The basic idea is that when $G$ has the Weyl property, its generators form a Coxeter system, meaning that the only relations between the generators are of the form $(\varphi_i\varphi_j)^{m_{ij}}=1$. As a consequence we can reduce the problem to dimension 2 and apply Theorem \ref{thm:1}. More generally this idea allows one to classify all models in arbitrary dimension for which the group is finite and the generators form a Coxeter system. We note that whenever $G$ is finite, Theorem \ref{thm:2} implies that $G$ is a Coxeter group, but the generators $\varphi_i$ do not always form a Coxeter system.

%We can go further and classify the models for which the group is finite and the generators form a Coxeter system, by reducing the problem to dimension two. Let us recall that $(\varphi_1, \ldots, \varphi_d)$ forms a Coxeter system. By definition, this means that the only relations between the generators are of the form $(\varphi_i\varphi_j)^{m_{ij}}=1$, where the integers $m_{ij}$ satisfy $m_{ii}=1$ and if $i\neq j$, $m_{ij}\geq 2$.

\section{Proof outline for Theorem~\ref{thm:1}}
\label{sec:3}

While the full proof would exceed the space available in this short note, we can sketch it as follows:
\begin{enumerate}
    \item\label{it1}We first establish a criterion for the finiteness of the group $\langle \varphi_1, \varphi_2 \rangle$ introduced in \eqref{eq:group_G_def}. To each model we associate a function $r:(0,1)\to(0,\infty)$, which is real-analytic in $t$ and expressed as the ratio of two elliptic integrals. The criterion can be stated as follows: the group is finite if and only if $r(t)=\frac{p}{q}$ is a rational constant; in that case, the order of the group is $2q$.
\item\label{it2}%Although the behavior of $r(t)$ is generally intricate (as it is a transcendental function), its limit as $t\to 0$ is easier to describe.
We show that $\lim_{t\to 0} r(t) = r_0$, where
    \begin{equation}
    \label{eq:possible_values_r0}
        r_0\in\left\{\frac{1}{2}, \frac{1}{3}, \frac{1}{4}, \frac{2}{3}, \frac{2}{5}, \frac{2}{7}, \frac{3}{4}, \frac{3}{5}, \frac{3}{7}, \frac{3}{8}, \frac{4}{7}, \frac{5}{7}, \frac{5}{8}\right\}.
    \end{equation}
    %Importantly, the value of $r_0$ depends only on the support of the step-set (that is, on the nonzero jumps), and not on the specific values of the weights.
    As a consequence, the only possible orders of the group are $4, 6, 8, 10, 14,$ and $16$. Surprisingly, $12$ does not appear as a possible order. 
    \item\label{it3}Using, on the one hand, that groups of orders $4$, $6$, and $8$ are already classified, and on the other hand, that the order of the group is at most $12$ (see \cite[Rem.~5.1]{HaSi-21}), the only remaining case is that of a group of order $10$. In this case we conclude by analysing the behaviour of $r(t)$ around $0$ that the only cases where the group has order 10 are those found in \cite{KaYa-15}. 
\end{enumerate}
We now give some more details on each item above. We start with \ref{it1}.
We use a detour via complex analysis. As it turns out, the group $\langle \varphi_1,\varphi_2\rangle$ in \eqref{eq:group_G_def} can not only be viewed as a group of birational transformations as in \eqref{eq:expression_generators}, but also as a group of automorphisms of the Riemann surface
\begin{equation*}
\{(x,y)\in(\mathbb C\cup{\infty})^2 : 1-t \chi_\mathcal S(x,y)=0\}.
\end{equation*}
For $t\in(0,1)$, and under hypothesis \ref{it:hypothesis_irreducible}, this surface has genus~$1$; it is therefore homeomorphic to a torus $\mathbb C/(\omega_1\mathbb Z + \omega_2\mathbb Z)$, where the fundamental periods $\omega_1,\omega_2$ of the lattice can be computed explicitly in terms of~$t$ as follows:
$\omega_1 = \alpha i K(k')$, $\omega_2 = \alpha K(k)$,
where $k$ is the elliptic modulus, which is an algebraic function of~$t$, $\alpha>0$ is a non-essential %constant
quantity (also algebraically dependent on $t$), and $K$ denotes the complete elliptic integral of the first kind; see, e.g., \cite[Sec.~6.3]{BeBMRa-21}. It is then shown that the group takes the form $\langle \omega\mapsto -\omega,\omega\mapsto -\omega+\omega_3\rangle$, with $\omega_3 = \alpha F(w,k)$, where $w$ is another algebraic function of~$t$ and $F(w,k)$ denotes the incomplete elliptic integral of the first kind. The function $r(t)$ is equal to the ratio $\frac{\omega_3}{\omega_2}=\frac{F(w,k)}{K(k)}$, and the group is finite if and only if $r(t)$ is a fixed rational number.

For point \ref{it2}, we proceed as follows: On the one hand, given the support of the step-set $\{(i,j):w(i,j)\neq 0\}$, we can directly write $w^2$ and $k^2$ as series in $\mathbb{C}[[\sqrt{t}]]$ with coefficients depending on the (non-zero) weighs $w(i,j)$. On the other hand, we can write $w^2$ and $k^2$ as series in the elliptic nome $q=\exp(i\pi\tau)=\exp(i\pi\frac{K(k)}{K'(k)})$ and $q^{r(t)}$, by inverting the relation $r=\frac{F(w,k)}{K(k)}$ using Jacobi elliptic functions: $w=\text{sn}(r K(k),k)=-i\text{sc}(r K(k),\sqrt{1-k^{2}})$, where the first equality is the classical inversion formula for $F(w,k)$, while the second is by the Jacobi transformation. Writing the Jacobi elliptic function $\text{sc}$ in terms of Jacobi theta functions yields the following expression for $w^2$:
%\aep{To check. My calculations led to the expressions below with $r/2$ replaced by $r$, but it seems correct as written, so maybe the definition of $r$ above should be multiplied by $2$.}
\begin{align*}
    w^{2}&=-\frac{\theta_{3}(q)^2\theta_{1}(\frac{r}{2}\tau,q)^2}{\theta_{4}(q)^2\theta_{2}(\frac{r}{2}\tau,q)^2}&&=\frac{(1+2q+2q^4+\cdots)^{2}(1-q^{r}-q^{2-r}+q^{2+2r}+\cdots)^{2}}{(1-2q-2q^4+\cdots)^{2}(1+q^{r}+q^{2-r}+q^{2+2r}+\cdots)^{2}},\\
    k^{2}&=\frac{\theta_{4}(q)^4}{\theta_{3}(q)^4}&&=\frac{(1-2q-2q^4+\cdots)^{4}}{(1+2q+2q^4+\cdots)^{4}}.
\end{align*}

Using these series in $q$, we have $\frac{\log(1-w^{2})}{\log(1-k^{2})}\to r_{0}$ as $q\to 0$, or equivalently as $t\to 0$. Then using the series in $t$, we can identify this limit, and show that it only depends on the support of the step-set.

Finally, for point \ref{it3}, we consider the situation where the group has order $10$. By symmetry arguments, we only need to consider the case where $r(t)=\frac{2}{5}$ and $w(1,1)=0$, while $w(1,-1),w(-1,-1),w(-1,1)=1$ and $w(1,0),w(0,1)\neq0$. Then from the series in $q$, we have
\[\frac{\left(385 \, w^{6} - 1415 \, w^{4} + 1835 \, w^2 - 869\right)}{\left(w^2 - 1\right)^{3}} + 32 \frac{\left(k^{2} - 1\right) \left(8 \, w^2 - 13\right)}{\left(w^2 - 1\right)^{6}} - 256 \frac{\left(k^2 - 1\right)^{2}}{\left(w^2 - 1\right)^{8}}=o(1),\]
as $q\to 0$. On the other hand, writing this expression directly as a series in $t$, the coefficient for $t^{-3}$, $t^{-2}$, $t^{-1}$ and $t^0$ are all rational functions of the square-roots of the unknown weights $w(1,0)$, $w(0,1)$, $w(-1,0)$, $w(0,-1)$. The only case where these coefficients are all $0$ is when $w(1,0)=w(0,1)=1$ and $w(-1,0)=w(0,-1)=2$.
This corresponds, after a central reflection, to the rightmost model in Figure~\ref{fig:step_sets}.

\section{Proof of Theorem~\ref{thm:2}}

It was shown in \cite{GoHuRa-25} that there is a surjective morphism $\phi: G\to H$, implying that whenever $G$ is finite, $H$ is also finite, and more precisely
\begin{equation*} 
%\label{GH}
H  \equiv \frac{G}{\ker\phi}.
\end{equation*}
Theorem \ref{thm:2} is then equivalent to the statement that $\ker\phi$ is trivial whenever $G$ is finite. Actually we prove the stronger statement that $\ker\phi$ is torsion-free, that is every non-identity element of $\ker\phi$ has infinite order.

We start by defining the transformation. From the definition \eqref{eq:chi_critical_point} of $x_{0}$, one can show that $x_{0}$ is fixed by each generator $\varphi_{j}$ of $G$, and hence by every element $g\in G$. This allows us to define the morphism
\begin{equation} \label{defJ}
J : \left| \begin{array}{ccc} 
G & \rightarrow & GL_d(\mathbb R) \\
g & \mapsto & \Jac_{\boldsymbol{x_0}} g
\end{array} \right.
\end{equation}
%with $\boldsymbol{x_0}$ defined in \eqref{eq:chi_critical_point}, and each $g\in G$ is a vector of rational functions of $x_{1},\ldots,x_{d}$. By a direct computation, one checks that for all $g \in G$, $g(\boldsymbol{x_0})=\boldsymbol{x_0}$, which implies that $J$ is a morphism.
In \cite{GoHuRa-25}, it is shown that $\Ima J$ and $H$ are isomorphic, which ensures the existence of a surjective morphism $\phi : G \to H$ for which $\ker J$ and $\ker\phi$ are isomorphic. Hence it suffices to show that $\ker J$ is torsion free.  %In particular,
%\begin{equation*} 
%\label{GH}
%H  \equiv \frac{G}{\ker\phi}.
%\end{equation*}
%Actually, we establish the following stronger version of Theorem~\ref{thm:2}:
%Let $g \in \ker\phi \setminus \{\Id\}$. Then the order of $g$ in $G$ is infinite. In particular, if $G$ is finite, then $G$ and $H$ are isomorphic.
This statement follows directly from the following lemma:
\begin{lem} 
\label{main_lem}
Let $\Om \subset \mathbb R^d$ be a connected open set and $f: \Om \to \Om$  be a $\mathcal C^2$-mapping such that 
\begin{itemize}
    \item there exists a positive integer $n$ such that $f^n=\Id$ (implying that $f$ is a $\mathcal C^2$-diffeomorphism);
    \item $f$ admits a fixed point $X \in \Om$ with $df_{X}=\Id$.
\end{itemize}
Then $f = \Id$.
\end{lem}
Indeed, let $g\in\ker J$. Then $g$ has the fixed point $x_{0}$ at which $dg_{x_{0}}$ is the identity. If additionally the order $n$ of $g$ is finite, then $g$ satisfies the assumptions of Lemma~\ref{main_lem} with $\Om= (0,\infty)^d$ and $X= \boldsymbol{x_0}$ which proves that $g=\Id$. Hence, the only element of $\ker J$ of finite order is the identity, so $\ker J$ is torsion free, as required. It remains to prove Lemma~\ref{main_lem}. Unlike Sections~\ref{sec:3} and \ref{sec:5}, we can present the complete proof here.
\begin{proof}
Let 
\begin{equation*}
    F:=\{x \in \Om| f(x)=x \; \hbox{ and } df_x = \Id\}.
\end{equation*} 
From the assumptions, we have $F \neq \emptyset$. Moreover, $F$ is clearly a closed set. The proof will be complete once we show that $F$ is also open, which, by the connectedness of $\Omega$, will imply that $F = \Omega$.
Hence it suffices to show that for $x \in F$ there is some $\gamma>0$ satisfying $B(x,\gamma)\subset F$. For the rest of the proof we fix $x\in F$. The idea is to show that for $z$ in an open ball centered at $x$, we have $f^{n}(z)-z\approx n(f(z)-z)$; then, since the left-hand side is $0$ by the first assumption of Lemma~\ref{main_lem}, we must have $f(z)=z$.

To make this approximation precise, we fix $\alpha>0$ such that $B(x,\alpha)\subset\Om$. We use the Taylor-Lagrange theorem, together with the fact that $f^k$ is of class $\mathcal C^2$, which implies that there exists $C_0>0$ such that for all $y,z \in B(x,\alpha)$ and all $k \in \{0,\ldots,n-1\}$ (where we recall that $n$ is the order of $f$ in the diffeomorphism group of $\Om$),
\begin{equation*}\|f^k(y) - f^k(z) -d(f^k)_z(y-z)\|\leq C_0 \|y-z\|^2.
\end{equation*}
%\notea_{(peut-être qu'il fallait l'argument que $\|d^2(f^k)_y\| \leq C$ aussi)}
In particular, if $f(z),z \in B(x,\alpha)$ and $v= f(z)-z$ then
\begin{equation}\label{T-Lfz}\|f^{k+1}(z) - f^k(z) -d(f^k)_z(v)\|\leq C_0 \|v\|^2.
\end{equation}
Now let $\ep\in(0,\alpha)$ be sufficiently small that $\varepsilon(1+2C_{0})<1$. Let $k \in \{0,\ldots,n-1\}$. Since $f^k$ is differentiable at $x$, $df_x=\Id$, and $x$ is a fixed point of $f$, it follows that $d(f^k)_x =\Id$. Since $f$ is of class $\mathcal C^1$, there exists $\de_k>0$ %and $C>0$
such that $B(x,\de_k) \subset \Om$ and such that for all $z \in B(x,\de_k)$,
\begin{equation} \label{normdf}
\|d(f^k)_z -\Id\| \leq \ep . %\hbox{ and }   \|d^2(f^k)_y\| \leq C.
\end{equation}
Finally, by continuity of $f^k$,  we can choose $\gamma_k\in(0,\de_k)$ such that 
\begin{equation*}
   f^k (B(x,\gamma_k)) \subset B(x,\min( \ep,\de_k)).   
\end{equation*}
Set $\gamma=\min_k (\gamma_k,\ep) $. We will show that $B(x,\gamma) \subset F$.
Combining \eqref{T-Lfz} and \eqref{normdf} applied at $v=f(z)-z$ yields %and Taylor-Lagrange's theorem, it holds that there exists $C_0>0$ such that for all  $v \in \R^d$, all $z \in B(x,\gamma)$. and all $k \in \{0,\cdots,n-1\}$,
%\begin{equation} \label{fz-z}
%\begin{aligned}
%\|f^{k}(z+v) -f^k(z) -v \| & \leq \|f^k(z+v) - f^k(z) -d(f^k)_z(v)\| + \|d(f^k)_z(v) -v\| \leq C_0 \|v\|^2 +\ep \|v\|. 
%\end{aligned}
%\end{equation}
\begin{equation} \label{fz-z}
\|f^{k+1}(z) -f^k(z) -v \| \leq C_0 \|v\|^2 +\ep \|v\|. 
\end{equation}
%Choose $k \in \{0, \cdots, n-1\}$, $z \in B(x,\gamma)$ and set $v=f(z) -z$. Since 
%$$f^{k+1} (z)-f^k(z)=f^k(z+v) - f^k(z),$$ 
By summing \eqref{fz-z} from $k=0$ to $k=n-1$, we obtain that, since $f^n=\Id$,
%$$f^{k+1} (z)-f^k(z) = v + v_k$$
%where $\|v_k\| \leq C_0 \|v\|^2  + \ep \|v\|$. 
%Summing for $k=0$ to $k=n$, we obtain 
\begin{equation} \label{0=v}
n\|v \|=\|f^{n}(z) -z -n v \| \leq\sum_{k=0}^{n-1}\|f^{k+1}(z) -f^k(z) -v \| \leq n(C_0 \|v\|^2 +\ep \|v\|).
\end{equation}
%where $\|w\| \leq n  (C_0 \|v\|^2  + \ep \|v\|$). 
Now observe that, by definition of $\gamma$, we have $f(z),z \in B(x,\ep)$, which implies that
$\| v\| \leq 2 \ep$, and thus the right-hand side of \eqref{0=v} is at most $(2 C_0 \ep + \ep) n\|v\|$.
But our choice of $\ep$ implies $(2 C_0 \ep + \ep)<1$, which makes \eqref{0=v} possible only if $v=0$, that is, $f(z)=z$. Hence $f(z)=z$ for all $z\in B(x,\gamma)$. This implies that $df_{z}=\Id$ for $z\in B(x,\gamma)$, and therefore $B(x,\gamma)\subset F$. This completes the proof of Lemma \ref{main_lem}.
\end{proof}

\section{Proof outline for Theorem~\ref{thm:3}}
\label{sec:5}

We assume in this section that $d=3$ and thus, $\chi_\mathcal S$ can be written as
\begin{equation*}
    \chi_\mathcal S (x,y,z) = \sum_{(i,j,k) \in \{-1,0,1\}^3} w(i,j,k)x^iy^jz^k,
\end{equation*}
where $(x,y,z)$ are the standard coordinates on $\mathbb{R}^3$. %A first observation is that, by Theorem~\ref{thm:2}, if $G$ is finite, then it must be isomorphic to a finite, 3 generator Coxeter group, that is, it must be one of the following groups: 
%\begin{equation*}
%\frac{\mathbb{Z}}{2\mathbb{Z}}\times\frac{\mathbb{Z}}{2\mathbb{Z}}\times\frac{\mathbb{Z}}{2\mathbb{Z}},\quad  \frac{\mathbb{Z}}{2\mathbb{Z}}\times D_{2k},\quad A_3, \quad B_3, \quad H_3, 
%\end{equation*} 
%where $D_{2k}$ is the dihedral group of order $2k$, and $A_3$, $B_3$, and $H_3$ are from the classification of finite Coxeter groups; see \cite[Chap.~2]{Hu-90}.
Denote by $m_{ij}$ the order of $\varphi_{i}\varphi_j$ in $G$. In this section, the main tool we use is the following result:
\begin{proposition}
\label{class_main}
The group $G$ is a finite group with the Weyl property if and only if
\begin{enumerate} 
\item \label{condition1} up to a permutation of the coordinates $(x,y,z)$, the triplet $(m_{12},m_{13},m_{23})$ belongs to the following list:
\begin{equation}
    \label{eq:list_triplets}
    \{(2,2,k),(3,2,3),(3,2,4),(3,2,5)| k \in \mathbb{N}^*\};
\end{equation}
 \item \label{condition2} for all $i \not=j$, the coefficient $a_{ij}$ of $\Delta$ in \eqref{eq:expression_covariance_matrix} is equal to $-\cos\left( \frac{\pi}{m_{ij}} \right)$.
 \end{enumerate}
\end{proposition}
%\aep{The equations in condition \ref{condition2} are wrong, they should be something like $\pi/m_{ij}=\arccos(a_{ij})$. in \cite[Prop.~4]{GoHuRa-25} it actually says $\pi/m_{ij}=-\arccos(a_{ij})$, but the sign doesn't seem to make sense.} \eh{I agree and made a change : $\cos$ instead of $\arccos$. What do you mean by "the sign doesn't seem to make sense " ? The sign allows to distinguish for instance $2\pi/3$ and $\pi /3$. Maybe I misunderstood what you mean}
The above result is a direct consequence of the results in \cite{GoHuRa-25} combined with Theorem \ref{thm:2} and the classification of finite Coxeter groups. Indeed, if conditions \ref{condition1} and \ref{condition2} both hold, then the classification of finite Coxeter groups (see \cite[Sec.~6.2]{GoHuRa-25}), combined with \cite[Prop.~4]{GoHuRa-25} imply that the group $H$ has the presentation
$$\left\{ r_1^2, r_2^2, r_3^2, (r_1 r_2)^{m_{12}},(r_1 r_3)^{m_{13}},(r_2 r_3)^{m_{23}} \right\},$$
where $r_{1},r_{2},r_{3}$ are the generators of $H$ defined in Section~\ref{reflection_group}, and that $T$ in \eqref{eq:def_new_domain} is a Weyl chamber of the group. Then condition \ref{condition1} implies that $G$ is finite and therefore isomorphic to $H$ by Theorem~\ref{thm:2}. In the other direction, if $G$ is finite and has the Weyl property then by Theorem \ref{thm:2}, $H$ is isomorphic to $G$ and $H$ is the reflection group associated to a Weyl chamber $T$, so it is a finite Coxeter system. Moreover, the order of $r_{i}r_{j}$ in $H$ is $m_{ij}$, which is the order of $\varphi_{i}\varphi_{j}$ in $G$. Then condition \ref{condition1} follows from the classification of finite Coxeter groups while condition \ref{condition2} follows from \cite[Prop.~4]{GoHuRa-25}.

%The above result is a direct consequence of the results in \cite{GoHuRa-25} combined with Theorem \ref{thm:2} and the classification of finite Coxeter groups. Indeed, the classification of finite Coxeter groups (see \cite[Sec.~6.2]{GoHuRa-25}) implies that condition \ref{condition1} is equivalent to the statement that $G$, along with its generators $\varphi_{1},\varphi_{2},\varphi_{3}$ is a finite Coxeter system. By theorem \ref{thm:2}, this is equivalent to the statement that $G$ is finite and  that  when $G$ is a finite group with the Weyl property, we know from \cite[Prop.~4]{GoHuRa-25} that $H$, along with its generators $r_{1},r_{2},r_{3}$ (defined in Section~\ref{reflection_group}) forms a Coxeter system. Then from the classification of finite Coxeter groups (see \cite[Sec.~6.2]{GoHuRa-25}), the group $H$ has the presentation
%$$\left\{ r_1^2, r_2^2, r_3^2, (r_1 r_2)^{m_{12}},(r_1 r_3)^{m_{13}},(r_2 r_3)^{m_{23}} \right\},$$
%where $(m_{12},m_{13},m_{23})$ are a permutation of one of the triples in \eqref{eq:list_triplets}. Moreover, having the Weyl property is equivalent to satisfying the second condition of the proposition. I the other direction, if the two conditions are satisfied then $H$, along with its generators $r_{1},r_{2},r_{3}$ forms the coxeter system above
%where we recall that $(r_1,r_2,r_3)$ are the generators of $H$ defined in Section~\ref{reflection_group}. 
%Then \cite[Prop.~11]{GoHuRa-25} implies that $G$ is finite (and therefore isomorphic to $H$ by Theorem~\ref{thm:2}).
Note that if $(m_{12},m_{13},m_{23})$ belongs to the list~\eqref{eq:list_triplets},
then $G$ is respectively isomorphic to $\frac{\mathbb{Z}}{2\mathbb{Z}}\times D_{2k}$, $A_3$, $B_3$, or $H_3$.

The next step is to use the classification of two-dimensional finite group models (according to Theorem \ref{thm:1} and Kauers and Yatchak’s paper~\cite{KaYa-15}) to determine when conditions \ref{condition1} and \ref{condition2} of Proposition~\ref{class_main} hold. % We now explain how this can be done.
To do this, we first prescribe the value of the triplet $(m_{12},m_{13},m_{23})$ (and thus $a_{ij} = -\cos\left( \frac{\pi}{m_{ij}} \right)$ ). Considering $z$ to be fixed, we define 
\begin{equation}
    \label{eq:def_chi-z}
    \chi_z (x,y):= \chi_\mathcal S(x,y,z) - ( w(0,0,1) z + w(0,0,-1) \overline{z}).
\end{equation}
Note that removing the constant term $w(0,0,1)z + w(0,0,-1)\overline{z}$ in~\eqref{eq:def_chi-z} does not affect the nature of the two-dimensional model in the variables $x,y$. Hence, for any $z>0$, the inventory $\chi_z (x,y)$ corresponds to a two-dimensional model with group of order $2m_{1,2}$, and coefficients %can be considered as a model in dimension $2$ and thus we denote by 
$$w^{xy}(i,j) :=  w(i,j,1) z +w(i,j,0) +  w(i,j,-1) \overline{z}.$$
Moreover, the value $a_{12}$ satisfies condition \ref{condition2}, which immediately allows us eliminate the possibility that $m_{12}=5$.  
%its  coefficients.
%since then   
%\begin{multline*}
%  \chi_z (x,y) =   w^{xy}(1,1)xy + w^{xy}(1,0) x + w^{xy}(1,-1) x \overline{y} + w^{xy}(0,1) y + w^{xy}(0,-1) \overline{y} \\  + w^{xy}(-1,1) \overline{x} y+w^{xy}(-1,0)\overline{x} + w^{xy}(-1,-1) \overline{x}\overline{y}.  
%\end{multline*}
We denote by $G_z$ the group associated with the model $\chi_z$.
Analogously, we define $\chi_x$, $\chi_y$, $G_x$, and $G_y$.
With this notation, and since in two dimensions the only finite combinatorial groups are $D_4$, $D_6$, $D_8$, and $D_{10}$, Proposition~\ref{class_main} can be restated as follows:
\begin{proposition} \label{class_main2}
The group $G$ is a finite group with the Weyl property if and only if
\begin{enumerate} 
\item up to a permutation of the coordinates $(x,y,z)$, the triple of groups $(G_z, G_y, G_x)$ appears in the following list:
$$\{(D_4,D_4,D_{2k}),(D_6,D_4,D_6),(D_6,D_4,D_8) : 2 \leq k  \leq 4 \};$$
 \item $a_{12} = -\cos\left( \frac{2\pi}{|G_z| }\right)$,  $a_{13} = -\cos\left( \frac{2\pi}{|G_y|} \right)$  and  $a_{23} = -\cos\left( \frac{2\pi}{|G_x|} \right)$.
 \end{enumerate}
\end{proposition}
So constructing a finite group $G$ with the Weyl property is equivalent to prescribing $G_x$, $G_y$ and $G_z$, which gives rise to equations on the coefficients $w^{xy}(i,j)$, $w^{xz}(i,j)$ and $w^{yz}(i,j)$ (from the classification in dimension $2$, see Section \ref{parag_dim_2}) which in turn lead to equations on the weights $w(i,j,k)$. Note that, to fulfil the second condition in Proposition~\ref{class_main2}, the classification in dimension~2 must only take into account the models whose covariance matrix has non-positive off-diagonal entries.
% (and just start to check which ones have this property). 

%Let us conclude this section by presenting a family of models obtained through this method, which is not a central weighting of any unweighted model. For this particular model, we have $G=A_{3}$ and more precisely $(m_{12}, m_{13}, m_{23})=(3,2,3)$. %(all classified using a computer, see \cite{KaWa-17}) by a transformation of the type
%$\chi_\mathcal S' = a \chi_\mathcal S(bx,cy,dz)$ (with $a,b,c,d>0$).

%\begin{proposition} \label{1aS2a'lem}
% Assume that 
% \begin{enumerate}
%     \item $w(0,0,-1) \geq0$ and $w(i,j,k) > 0$ if  $(i,j,k) $ belongs to
%     $$ \{(0,-1,0),(-1,1,-1),(-1,1,0),(0,-1,-1),(0,-1,1),(1,0,-1),(1,0,0)\}, $$
%%\item   $w(0,0,-1) \geq0$; 
%\item $w(i,j,k) = 0$ otherwise; 
%\item finally $w(0,-1,0) w(-1,1,-1)=2w(-1,1,0)w(0,-1,-1)$,   $w(-1,1,0)w(0,-1,0)=$ $2w(-1,1,-1)w(0,-1,1)$,      and $w(1,0,-1)w(-1,1,0)=w(1,0,0)w(-1,1,-1)$. 
% \end{enumerate}
% Then $G$ is finite, isomorphic to $A_3$ and has the Weyl property.
% \end{proposition}

Let us conclude this section by giving the classification for the Coxeter groups $A_3$ and $B_3$.
For any inventory $\chi_\mathcal S(x,y,z)$ corresponding to a finite group $G$ with the Weyl property, inventories defined by permuting $x,y,z$, as well as the inventories $a \chi_\mathcal S(bx,cy,dz)$ and $\chi_\mathcal S(\overline x,\overline y,\overline z)$ also have these properties. Up to these transformations, we show that there are exactly two families of models for which $(m_{12}, m_{13}, m_{23})=(3,2,3)$, so $G=A_{3}$. The first is represented by
$$\chi_\mathcal S(x,y,z)=\frac{1}{yz}+\frac{2}{y}+\frac{z}{y}+\frac{x}{z}+x+\frac{y}{xz}+\frac{y}{x}+\frac{c}{z},$$
for $c\geq 0$, while the second is represented by
\begin{align*}\chi_\mathcal S(x,y,z)&=a\left(x+\frac{y}{x}+\frac{z}{y}+\frac{1}{z}\right)+b\left(\frac{1}{x}+\frac{x}{y}+\frac{y}{z}+z\right)+c\left(y+\frac{1}{y}+\frac{xz}{y}+\frac{y}{xz}+\frac{z}{x}+\frac{x}{z}\right),\end{align*}
where $a,b,c\geq 0$ and $a,b,c$ are not all $0$.
We have also fully classified the models with $(m_{12}, m_{13}, m_{23})=(3,2,4)$, so the group $G=B_{3}$. In this case, all models are related by the symmetries described above to one of two unweighted models:
\begin{equation*}\chi_\mathcal S(x,y,z)=\frac{y}{z}+\frac{z}{y}+\frac{x}{y}+\frac{y}{x}+\frac{1}{x}+x,\end{equation*}
and
\begin{equation*}\chi_\mathcal S(x,y,z)=\frac{x}{z}+\frac{z}{x}+\frac{y}{z}+\frac{z}{y}+\frac{y}{xz}+\frac{xz}{y}+z+\frac{1}{z}.\end{equation*}

\vspace{-6mm}

%\acknowledgements{}

%% if you use biblatex then this generates the bibliography
%% if you use some other method then remove this and do it your own way
\printbibliography

\end{document}